\catcode`\@=11
\def\citen#1{\if@filesw \immediate\write \@auxout {\string\citation{#1}}\fi%
\@tempcntb\m@ne \let\@h@ld\relax \def\@citea{}%
\@for \@citeb:=#1\do {\@ifundefined {b@\@citeb}%
    {\@h@ld\@citea\@tempcntb\m@ne{\bf ?}%
    \@warning {Citation `\@citeb ' on page \thepage \space undefined}}%
    {\@tempcnta\@tempcntb \advance\@tempcnta\@ne
    \setbox\z@\hbox\bgroup\ifcat0\csname b@\@citeb \endcsname \relax
    \egroup \@tempcntb\number\csname b@\@citeb \endcsname \relax
    \else \egroup \@tempcntb\m@ne \fi \ifnum\@tempcnta=\@tempcntb
    \ifx\@h@ld\relax \edef \@h@ld{\@citea\csname b@\@citeb\endcsname}%
    \else \edef\@h@ld{\hbox{--}\penalty\@highpenalty
    \csname b@\@citeb\endcsname}\fi
    \else \@h@ld\@citea\csname b@\@citeb \endcsname \let\@h@ld\relax \fi}%
\def\@citea{,\penalty\@highpenalty\hskip.13em plus.13em minus.13em}}\@h@ld}
\def\@citex[#1]#2{\@cite{\citen{#2}}{#1}}%
\def\@cite#1#2{\leavevmode\unskip\ifnum\lastpenalty=\z@\penalty\@highpenalty\fi%
  \ [{\multiply\@highpenalty 3 #1%
  \if@tempswa,\penalty\@highpenalty\ #2\fi}]}   %
\makeatother 

\catcode`\@=12

\def\calg  {{\mathcal G}}
\def\calh  {{\mathcal H}}
\def\call  {{\mathcal L}}
\def\calm  {{\mathcal M}}
\def\caln  {{\mathcal N}}
\def\calr  {{\mathcal R}}
\def\cals  {{\mathcal S}}
\def\calt  {{\mathcal T}}
\def\calv  {{\mathcal V}}
\def\cft   {conformal field theory}
\def\cfts  {conformal field theories}
\def\complex {{\mathbb C}}
\def\g       {{\mathfrak g}}
\def\vvir    {v_{{\rm vir}}}
\def\zet     {{\mathbb Z}}

\documentclass[12pt]{article} \usepackage{amssymb,amsfonts,latexsym}

\begin{document}


\begin{flushright}  {~} \\[-1cm]
{\sf math/0004034}\\{\sf PAR-LPTHE 00-09}\\
{\sf April 2000} \end{flushright}

\begin{center} \vskip 14mm
{\Large\bf D-BRANE CONFORMAL FIELD THEORY}\\[4mm]
{\Large\bf AND BUNDLES OF CONFORMAL BLOCKS}\\[20mm]
{\large J\"urgen Fuchs$\;^1$, \ and
Christoph Schweigert$\;^2$}
\\[6mm]
$^1\;$ Institutionen f\"or fysik\\
Universitetsgatan 1\\ S\,--\,651\,88\, Karlstad\\[3.5mm]
$^2\;$ LPTHE, Universit\'e Paris VI\\
4 place Jussieu\\ F\,--\,75\,252\, Paris\, Cedex 05\\[3.5mm]
\end{center}
\vskip 18mm
\begin{quote}{\bf Abstract}\\[1mm]
Conformal blocks form a system of vector bundles over the moduli space of 
complex curves with marked points. We discuss various aspects of these 
bundles. In particular, we present conjectures about the dimensions
of sub-bundles. They imply a Verlinde formula for non-simply connected 
groups like $PGL(n,C)$.

We then explain how conformal blocks enter in the construction of conformal
field theories on surfaces with boundaries. Such surfaces naturally appear
in the conformal field theory description of string propagation in the
background of a D-brane. In this context, the sub-bundle structure of the
conformal blocks controls the structure of symmetry breaking boundary
conditions.
\end{quote}

 \newpage


\section{Introduction}

Two-dimensional conformal field theory plays a fundamental role in the theory
of two-dimensional critical systems of classical statistical mechanics,
in quasi one-dimensional condensed matter physics and in string theory.
Moreover, this field of research has repeatedly contributed substantially
to the interaction between mathematics and physics. The study of defects 
in systems of condensed matter physics, of percolation probabilities, 
and of (open) string perturbation theory in the background of certain string
solitons, the so-called D-branes, has recently driven theoretical
physicists to analyze conformal field theories on surfaces that may
have boundaries and\,/\,or can be non-orientable.

In the present contribution, we would like to describe some mathematical
aspects of this recent development and state some conjectures 
about the sub-bun\-dle structure of the bundles of conformal blocks.
In physics, this structure enters in the description of symmetry
breaking boundary conditions. In the special case of WZW models
our conjectures imply a Verlinde formula for non-simply
connected groups like $PGL(n,C)$.

\medskip

We are grateful to the organizers of the third European Congress of
Mathematics for the invitation to present these results.

\section{Chiral and full conformal field theory}

The investigation of boundary conditions makes it clear that it is
indispensable to distinguish carefully between {\em chiral} conformal field 
theory, CCFT, and {\em full} conformal field theory, CFT. Although a CFT is
constructed from a CCFT, both types of theories are rather different
in nature. Mathematicians typically have in 
mind {\em chiral\/} \cft, CCFT, when they use the term \cft; one goal of the 
present note is to draw their attention to this difference. In physics, both
CCFT and CFT have found applications: CFT in string theory and the study
of two-dimensional critical systems, CCFT in the description of 
quantum Hall fluids (for a recent review see \cite{fpsw}).

A CCFT is defined on a closed, compact complex curve $\hat X$; a CFT, in 
contrast, is defined on a real two-dimensional manifold $X$ that is endowed 
with a conformal structure, i.e.\ a class of metrics modulo local rescalings.
While $\hat X$ is oriented and has empty boundary, $X$ is allowed to
have a boundary and/or to be non-orientable.

To any such $X$, one finds a complex curve $\hat X$,
the {\em  double\/}, together with an anti-conformal involution $\sigma$ of
$\hat X$  such that $X$ is isomorphic to the quotient of $\hat X$
by $\sigma$. (When $X$ has empty boundary, $\hat X$ is just the total
space of the orientation bundle of $X$.) $\hat X$ is not 
necessarily connected; also, $\sigma$ does not necessarily
act freely: fixed points of $\sigma$ correspond to boundary points of $\hat X$.

A double can also be constructed when $X$ has dimension higher than two,
and it will be an oriented conformal manifold. However, only in two dimensions
a conformal structure plus an orientation are equivalent to a holomorphic
structure; in higher dimensions, additional integrability constraints have
to be fulfilled. This simple fact implies that specifically in two dimensions
the powerful tools of the theory of holomorphic functions are at one's disposal. 

We will start with a discussion of chiral conformal field theory, exhibiting 
in particular some aspects of the bundles of conformal blocks. For instance,  
we will explain conjectures about the sub-bundle structure of these bundles. 
Full CFT will be discussed in Section \ref{s.3}. The general idea is to relate 
full CFT on $X$ to chiral conformal field theory on the double $\hat X$. 

\section {Chiral conformal field theory}

\subsection{Conformal vertex operator algebras}

The definition of a chiral conformal field theory starts with a 
conformal vertex operator algebra (VOA) $(\calh_\Omega,v_\Omega,Y,\vvir)$. 
Here $\calh_\Omega$ is a complex vector space, $v_\Omega$ and $\vvir$ are 
vectors in $\calh_\Omega$, and 
\begin{equation}
Y: \quad \calh_\Omega \to {\rm End}(\calh_\Omega) [[z,z^{-1}]]
\end{equation}
is a map into the space of Laurent series in a formal variable $z$ with
values in ${\rm End}(\calh_\Omega)$, called a {\em field-state correspondence\/}.
One imposes a weakened version of commutativity, called
{\em locality\/}: for each pair of vectors $v,w\,{\in}\,\calh_\Omega$, there
is an integer $N\,{=}\,N(v,w)$ such that
\begin{equation}
(z-w)^N Y(v,z) Y(w,b) =
(z-w)^N Y(w,b) Y(v,z) 
\label{1}
\end{equation}
in the sense of formal power series. For precise definitions and an introduction
to VOAs, see e.g.\ \cite{KAc4}. 

One way to think of the field-state correspondence
is as a family of products on the vector space $\calh_\Omega$ 
that depend on the formal variable $z$:
\begin{equation}
v \star_z w = Y(v,z) w  \,. 
\end{equation}
It is remarkable that ``commutativity'' (\ref{1}) implies a kind
of associativity: one has
\begin{equation}
a\star_z (b\star_w c) = (a\star_{z-w} b)\star_w c \,.
\end{equation}
for all $a,b,c\,{\in}\,\calh_\Omega$.
The vector $\Omega$ is called the vacuum vector. Under field-state 
correspondence, it is mapped to the identity map, 
$Y(\Omega,z)={\rm id}$.
The Virasoro vector $\vvir$ has the property that the modes $L_n$ in the
expansion of the corresponding field, the {\em stress energy tensor\/},
\begin{equation}
Y(\vvir,z) = \sum_{n\in\zet} L_n\, z^{-n-2} 
\end{equation}
span an infinite-dimensional Lie algebra that is isomorphic to the Virasoro
algebra.

At first sight a VOA might seem a rather ad hoc kind
of algebraic structure. We would like to stress, however, that it is a very
natural concept indeed. Apart from their motivation from physics
(for a review see \cite{KAc4}), VOAs can be naturally characterized in
a category theoretic framework as singular commutative rings in certain
categories \cite{borc}. Examples of VOAs can be constructed from various
in\-fini\-te-di\-mensional Lie algebras like the Heisenberg algebra,
the Virasoro algebra or untwisted affine Lie algebras; the latter case gives
rise to the so-called WZW models. Moreover, the chiral de Rham complex 
\cite{masv} allows to associate to any complex variety a VOA as an invariant. 

For every algebraic structure, it is natural to investigate its representation
theory. Thus denote by $I$ the set of equivalence classes of irreducible 
representations of a given VOA, and $\calh_\mu$ a representative for $\mu\,{\in}
\,I$. Elements of $I$ will also be called labels. Any irreducible representation
carries in particular a representation of the Virasoro algebra. Usually,
one requires that $L_0$ acts by semi-simple operators.\,%
\footnote{~This requirement is weakened in so-called logarithmic \cfts.}
The axioms of a VOA then imply that the spectrum of $L_0$ is integrally spaced
and bounded from below. The lowest eigenvalue of $L_0$ in $\calh_\mu$ is called
the {\em conformal weight\/} $\Delta_\mu$.

A special representation is given by the vector space $\calh_\Omega$ of the
VOA itself; it is called the vacuum representation, and the corresponding 
label in $I$ will be denoted by $\Omega$.
A VOA is called rational iff the set $I$ is finite and every representation
is completely reducible. The representation theory of VOAs is a model
dependent problem that can be quite intricate; see e.g.\ \cite{fusS4,xu} for 
the so-called coset \cfts.

\subsection{Conformal blocks}

VOAs formalize aspects of the `observables' of chiral two-dimensional 
field theories. As a consequence, the algebraic structure of a VOA fits 
together very well with certain two-dimensional geometries -- more precisely,
since the theory is chiral, with the geometry of complex curves. 
To see how this works out, we fix a complex 
curve $\hat X\,{=}\,\hat X_g$ of genus $g$, and $m$ distinct smooth 
points $\vec z\,{=}\,(z_1,...\,,z_m)$ on $\hat X$, the insertion points.
A VOA should be thought of as being associated to a (formal) punctured
disc with coordinate $z$; to the data $\hat X, \vec z$ one can associate
a global version of the VOA.\,%
\footnote{~This algebra has been named chiral algebra in \cite{beDr3}; 
unfortunately this
does not agree with terminology in physics, where a chiral algebra is a VOA.}
In the case of $\hat X\,{=}\,\complex P^1$ with three marked points, this global
version also leads to the coproduct for the VOA that was defined in \cite{fefk3}.

When we fix additional data, labels $\vec\lambda\,{=}\,\lambda_1,\lambda_2,
...\,,\lambda_m$ and local coordinates $\xi_1,\xi_2,...\,,\xi_m$ for each 
insertion point, this global algebra acts on the algebraic dual
\begin{equation}
(\calh_{\lambda_1}\otimes_{\complex}^{}\cdots\otimes_{\complex}^{}
\calh_{\lambda_m})^*_{}  \,,
\end{equation}
of the tensor product of the modules $\calh_{\lambda_i}$. 
The vector space $V_{\vec\lambda}(\hat X_g)$ of conformal blocks 
is by definition the invariant subspace under this 
action.  We call a complex curve with the additional data we just described an
{\em extended curve\/}.

For the vertex operator algebra of the WZW model based on an untwisted affine 
Lie algebra $\g=\bar\g^{\scriptscriptstyle (1)}$ at level 
$k\,{\in}\,\zet_{\ge0}$, the vector spaces
$V_{\vec\lambda}(\hat X_g)$ also appear naturally in algebraic geometry.
The space $V_{\Omega}(\hat X_g)$, e.g., can be identified with the space of 
holomorphic sections in the $k$-th power of a line bundle $\call$ over the
moduli space of holomorphic $G$-connections over $\hat X$, where $G$ is the
simple, connected and simply connected complex Lie group with Lie algebra
$\bar\g$.

One can show that the vector spaces $V_{\vec\lambda}(\hat X_g)$ with fixed 
$\vec\lambda$ and $g$ fit together into a vector bundle 
$\calv_{\vec\lambda,g}$ over the moduli space $\calm_{g,m}$ of complex
curves of genus $g$ with $m$ marked points. There is no reason why this
vector bundle $\calv_{\vec\lambda,g}$ 
should be irreducible. Later, we will describe some sub-bundles 
which enter in the description of symmetry breaking boundary conditions.

The Virasoro element $\vvir$ can also be used to endow
this bundle with a pro\-jectively flat connection, the Knizhnik-Zamolodchikov
connection. As a result, the dimension of the vector space 
$V_{\vec\lambda}(\hat X_g)$ depends only on the genus $g$ and on $\vec\lambda$.
The Knizhnik-Zamolodchikov connection implies a projective action of the 
fundamental group of $\calm_{g,m}$, the mapping class group, on 
$V_{\vec\lambda}(\hat X_g)$. In the particular case of $g\,{=}\,1$ and with
one insertion of the vacuum $\Omega$, we obtain a 
projective representation of the modular group $PSL(2,\zet)$ on a complex
vector space of dimension $|I|$. In a natural basis, the
generator $T{:}\;\tau\,{\mapsto}\,\tau{+}1$ of $PSL(2,\zet)$
is represented by a unitary diagonal matrix $T_{\lambda,\mu}$, and 
$S{:}\;\tau\,{\mapsto}\,{-}1/\tau$ by a unitary symmetric matrix 
$S_{\lambda,\mu}$.

Up to this point, the structures we discussed are mathematically essentially
under control (see e.g.\ \cite{frbz}). There is, however, one further
crucial aspect of
conformal blocks, known as factorization. The curve $\hat X$ is allowed
to have ordinary double points. Blowing up such a point $p$ yields
a new curve $\hat X'$ with a projection to $\hat X$ under which $p$
has two pre-images $p'_\pm$. By factorization one means the existence
of a canonical isomorphism 
\begin{equation}
V_{\vec\lambda}(\hat X) \,\cong\, \bigoplus_{\mu\in I} 
V_{\vec\lambda\cup\{\mu,\mu^+\}}(\hat X')
\end{equation}
between the blocks on $\hat X$ and $\hat X'$.
This structure tightly links the system of bundles $\calv_{\vec\lambda,g}$
over the moduli spaces $\calm_{g,m}$ for different values of $g$ and $m$. 
One consequence is that one can express the rank of $\calv_{\vec\lambda,g}$ 
for all values of $m$ and $g$ in terms of the matrix $S$ that we encountered 
in the description of the action of the modular 
group. This results in the famous {\em Verlinde formula\/}, which reads
\begin{equation}  
{\rm rank}\, \calv_{\vec\lambda,g} = \sum_{\mu\in I}
\prod_{i=1}^m \, \frac{S_{\lambda_i, \mu}}{S_{\Omega,\mu}} \,\,
|S_{\Omega,\mu}|^{2-2g} \, . 
\label{3} \end{equation}
The matrix $S$ can be computed explicitly for concrete models; in the
case of WZW models, $S$ is given by the Kac-Peterson formula. The
combination of the Kac-Peterson formula for $S$ with the general
Verlinde formula (\ref3) then gives the Verlinde formula in the sense
of algebraic geometry \cite{beau,falt}.

\subsection{Sub-bundles}

We now discuss sub-bundles of the bundles $\calv_{\vec\lambda,g}$.
To this end, we introduce the {\em fusion rules\/} 
$\caln_{\lambda_1\lambda_2\lambda_3}$, which are dimensions of the three-point 
blocks at genus 0:
\begin{equation}
\caln_{\lambda_1,\lambda_2,\lambda_3} = {\rm rank}\,
\calv_{\lambda_1,\lambda_2,\lambda_3} \,.
\end{equation}
It turns out that $\caln_{\lambda_1,\lambda_2,\Omega}$ describes a permutation
of order two on $I$ which we denote by $\mu\,{\mapsto}\,\mu^+$. Consider the free 
$\zet$-module $\calr=\zet^{|I|}$ with a basis $\{\Phi_\mu\}$ labelled by 
$\mu\,{\in}\,I$. The definition
\begin{equation}
\Phi_{\mu_1} {\star}\,\Phi_{\mu_2} = \sum_{\mu_3\in I}
\caln_{\mu_1^{},\mu_2^{},\mu^+_3}\, \Phi_{\mu_3}
\end{equation}
turns $\calr$ into a ring, the {\em fusion ring\/}. 
Factorization and Knizhnik-Zamolodchikov connection imply that this 
ring is commutative, associative and semi-simple. Notice that by construction 
a fusion ring comes with a distinguished basis labelled by $I$.

Invertible elements of $\calr$ that are elements of the distinguished basis
are called simple currents. In the WZW case,
simple currents correspond (with one exception, appearing for $E_8$ at level 2),
to elements of the center of the Lie group $G$ \cite{jf15}. The group $\calg$ of
simple currents acts on $\calr$, and also on $I$. This action is in general
not free, and to each element $\mu\,{\in}\,I$ we associate its stabilizer
$\cals_\mu:= \{ \Phi_J\,{\in}\, \calg \,|\, \Phi_J{\star} \Phi_\mu\,{=}\,\Phi_\mu\}$.
Given an $m$-tuple of labels $\vec\lambda$, consider the subgroup 
$\cals_{\vec\lambda}^1$ of $\cals_{\lambda_1}{\times}\cdots{\times}\,
\cals_{\lambda_m}$ that consists of those elements $(J_1,J_2,...\,,J_m)$ 
whose product equals the unit element of the fusion ring, $\prod_i J_i= \Omega$.

Henceforth we will for simplicity assume that $H^2(\calg,\complex^{\,*}
_{{\phantom|}})\,{=}\,0$; for results in more general situations see 
\cite{fusS6,fuSc8}. One expects that then the group $\cals_{\vec\lambda}^1$ 
acts by bundle 
automorphisms on $\calv_{\vec\lambda,g}$; for WZW models and $g\,{=}\,0$, 
see \cite{fuSc8}. The invariant subspaces under this action form sub-bundles;
once an action of $\cals_{\vec\lambda}^1$ is fixed, they  are characterized
by their eigenvalues, i.e.\ by a character of $\cals_{\vec\lambda}^1$. 

A conjecture for the ranks of these sub-bundles was presented in 
\cite{fusS6,fuSc8}. It is actually more convenient to give the Fourier
transforms of these ranks which can be interpreted as traces of the action
of $\cals_{\vec\lambda}^1$ on the conformal blocks. For the trace of 
$(J_1,J_2,...,J_m)\,{\in}\,\cals_{\vec\lambda}^1$ we obtain
\begin{equation}  
{\rm Tr}_{V_{\vec\lambda,g}}^{} \Theta_{J_1,\ldots, J_m} = \sum_{\mu\in I} \,
\prod_{i=1}^m \,\, \frac{S^{J_i}_{\lambda_i, \mu}}{S_{\Omega,\mu}} \,
|S_{\Omega,\mu}|^{2-2g}_{} \, . 
\label{4} \end{equation}
Here $S^J$ is the matrix that describes the transformation under the modular
group of the one-point blocks on an elliptic curve with insertion $J$.
In particular, for $J\,{=}\,\Omega$, $S^J$ equals the ordinary matrix $S$
and we recover the Verlinde formula (\ref3) for $\Theta\,{=}\,{\rm id}$.

In the case of WZW models, there is a conjecture \cite{fusS6} for the explicit 
form of the matrix $S^J$: it is given by the S-matrix of Kac-Peterson form
for another Lie algebra, the so-called orbit Lie algebra \cite{fusS4,furs} 
associated to the data $(\g,J)$.\,%
\footnote{~Orbit Lie algebras also enter in the study
of moduli spaces of flat connections with non-simply connected structure
group on elliptic curves \cite{schW3}, and in the related problem of
determining almost commuting elements in Lie groups \cite{bofm}.} 
With this conjecture, formula (\ref4) is as concrete as the Verlinde formula
and can be implemented in a computer program, see
{\tt http:/$\!$/www.nikhef.nl/\~{}t58/kac.html}. In the case of WZW models,
these sub-bundles also enter in a Verlinde formula, in the sense of algebraic
geometry, for non-simply connected groups like $PGL(n,\complex)$, see
\cite{beau3}.

\section{Full conformal field theory}\label{s.3}

Having discussed some structures of CCFT, our next goal
is to construct full CFT on a real two-dimensional manifold $X$
with conformal structure, using CCFT on the double $\hat X$ of $X$. 
To this end, we have to prescribe some data for each insertion point, which
can lie either in the interior of $M$ or on the boundary. For each point
on the boundary of $M$, we choose a label and a local coordinate
$\xi$ such that the lift to the double $\hat X$ gives a local holomorphic
coordinate. We also choose an orientation for each component of the boundary.

For each point in the interior of $M$, a so-called bulk point, we choose 
a local coordinate such that on the double we obtain local 
holomorphic coordinates around the two pre-images. Notice that possible
choices of such coordinates come in two disjoint sets which are related by
a complex conjugation of the coordinate. We call them a local orientation of
the coordinate.

To specify the labels of the bulk fields, we need to fix one more datum:
an automorphism $\dot\omega$ of the fusion ring that preserves conformal
weight modulo integers, $T_{\dot\omega\lambda}\,{=}\,T_{\lambda}$.
The labels for bulk fields are now defined as equivalence classes 
of triples $(\lambda,or',or)$, where $\lambda$ is a label of the chiral CCFT, 
$or'$ is a local orientation of the coordinate, and $or$ is a local orientation 
of $X$ around the insertion point. The triples $(\lambda,or',or)$
and $(\dot\omega\lambda,-or',-or)$ with opposite orientations are
identified. We call a surface with this structure a {\em labelled surface\/}
(cf.\ \cite{fffs3}, where the notion was introduced in the category of
topological manifolds for the case when $\dot\omega$ is the conjugation).

Some thought shows that these data precisely allow to give the double, with
the pre-images of the insertion points, the structure of an extended 
curve, so that we can define conformal blocks. While the moduli space
$\calm_{X}$ of a labelled surface cannot be embedded into the moduli space 
$\calm_{\hat X}$ of its double as an extended curve, an embedding 
of the corresponding Teichm\"uller spaces $\calt$ does exist.
This is sufficient to give us conformal blocks as vector bundles 
$\calv_{\vec\lambda,g}$ over $\calm_X$. The involution $\sigma$ on $\hat X$ 
induces an involution $\sigma_*$ on $\calt_{\hat X}$. The subgroup of the 
mapping class group of $\hat X$ that commutes with $\sigma_*$ is called 
relative modular group \cite{bisa2}.
(In the case of a closed and oriented surface, $\hat X$ has two connected
components of opposite orientation each of which is isomorphic to $X$. The 
relative modular group is then isomorphic to the mapping class group of $X$.)

The central problem in the construction of a CFT from a CCFT is to specify,
for all choices of $\vec\lambda$ and all genera $g$,
{\em correlation functions\/} of the CFT as sections in 
$\calv_{\vec\lambda,g}$. For fixed labels $\vec\lambda$ and fixed $g$, 
these sections depend
on the positions of the insertion points and on the moduli of the conformal
structure. We impose the locality requirement: these sections have to be
genuine functions (rather than multivalued sections) of these parameters.
This includes in particular invariance under the relative modular group.
Moreover, we require compatibility with factorization (for a formulation 
of the latter requirement in 
the category of topological manifolds, see \cite{fffs3}).

Due to factorization, the calculation of any correlation function can
be reduced to four basic cases: the sphere with three points, the disc
with three boundary points, the disc with one boundary point and one interior
point, and the real projective plane with one point. Except for the last case,
they involve three marked points on each connected component of the double. 

Constraints on these particular correlators follow by considering situations
where four marked points appear on each connected component of $\hat X$.
In particular, from the consideration of four points on the boundary of the 
disc one can deduce that the correlator of three boundary fields
is proportional to the so-called fusing matrix \cite{bppz2,fffs}.
We stress that this is a model independent statement. Indeed, the
construction of a CFT from a CCFT can be formulated in a completely
model independent way; see \cite{fffs2,fffs3}, where a general prescription 
for the correlators has been given for the case when $\dot\omega$ is charge
conjugation.

Similarly, one can consider two bulk fields on a disc to derive
constraints on the correlator of one bulk field and one boundary field on a
disc with boundary condition $a$. To be more explicit, introduce a basis 
$\{ e^{(\alpha)}_{\lambda,\mu,\nu} \}$ for the space $V_{\lambda,\mu,\nu}$
of three-point blocks on the sphere; we are looking for a linear combination
\begin{equation}
\sum_{\alpha=1}^{\caln_{\lambda,\dot\omega\lambda,\mu}}
{}^{\alpha\!}R^a_{\lambda,\dot\omega\lambda,\mu}\,
e_{\lambda,\dot\omega\lambda,\mu}^{(\alpha)}
\end{equation}
that is compatible with factorization. It turns out that, with suitable 
normalizations, the factorization constraints imply that
\begin{equation}
R^a_{\lambda_1,\dot\omega\lambda_1,\Omega}
R^a_{\lambda_2,\dot\omega\lambda_2,\Omega} =
\sum_{\lambda_3} \tilde\caln_{\lambda_1,\lambda_2}^{\lambda_3}
R^a_{\lambda_3,\dot\omega\lambda_3,\Omega} \,.
\end{equation} 
This equation shows that the {\em reflection coefficients\/}
$R^a_{\lambda,\dot\omega\lambda,\Omega}$ form one-di\-men\-si\-onal 
representations of a certain algebra, called the {\em classifying algebra\/} 
\cite{fuSc5}. The problem of classifying boundary conditions 
is therefore reduced to the representation theory of the classifying algebra.

It turns out that, in the case when $\dot\omega$ is charge conjugation,
the structure constants $\tilde\caln_{\lambda_1,\lambda_2}^{\lambda_3}$
are just the fusion rules, i.e.\ the dimensions of the spaces of three-point
blocks. This result has been generalized both to more general choices of 
the fusion rule automorphism $\dot\omega$ \cite{fuSc5} and to the case
when the boundary conditions break some of the symmetries of the theory
\cite{fuSc10,fuSc11,fuSc12}. 
The structure constants of the classifying algebra can be expressed in terms 
of the modular transformation matrices $S^J$ that appear in (\ref4) and 
of characters of simple current groups \cite{fuSc11,fuSc12}. The 
resulting expressions can be shown to coincide with traces of the form (\ref4) 
for appropriate twisted intertwiner maps on the spaces of three-point blocks. 
Via this relationship, the sub-bundle structure of the bundles of conformal 
blocks controls boundary conditions in these more general cases, in particular 
boundary conditions that break some of the bulk symmetries. Solitonic 
solutions in field and string theory typically do not respect all symmetries, 
and D-brane solutions are no exception to this. The sub-bundle structure of 
conformal blocks is therefore an essential ingredient for the CFT approach to
D-brane physics.

 \newcommand\wb{\,\linebreak[0]} \def\wB {$\,$\wb}
 \newcommand\Bi[1]    {\bibitem{#1}}
 \newcommand\J[5]   {{\sl #5}, {#1} {#2} ({#3}) {#4} }
 \newcommand\Prep[2]  {{\sl #2}, preprint {#1}}
 \newcommand\BOOK[4]  {{\em #1\/} ({#2}, {#3} {#4})}
 \newcommand\inBO[7]  {in:\ {\em #1}, {#2}\ ({#3}, {#4} {#5}), p.\ {#6}}
 \def\atmp  {Adv.\wb Theor.\wb Math.\wb Phys.}
 \def\comp  {Com\-mun.\wb Math.\wb Phys.}
 \def\joag  {J.\wB Al\-ge\-bra\-ic\wB Geom.}
 \def\joal  {J.\wB Al\-ge\-bra}
 \def\kma     {Kac-Moody algebra}
 \def\maan  {Math.\wb Annal.}
 \newcommand\mbop[2] {\inBO{The Mathematical Beauty of Physics}
            {J.M.\ Drouffe and J.-B.\ Zuber, eds.} \WS\Si{1997} {{#1}}{{#2}} }
 \def\mpla  {Mod.\wb Phys.\wb Lett.\ A}
 \def\nuci  {Nuovo\wB Cim.}
 \def\nupb  {Nucl.\wb Phys.\ B}
 \def\phlb  {Phys.\wb Lett.\ B}
 \def\phrl  {Phys.\wb Rev.\wb Lett.}
 \def\phrp  {Phys.\wb Rep.}
 \def\WS     {{World Scientific}}
 \def\Si     {{Singapore}}
 \def\AMS    {{American Mathematical Society}}
 \def\PR     {{Providence}}

\end{document}